\documentclass[12pt,reqno]{amsart}
\usepackage{amsmath,amssymb, amscd, eufrak}

\newcommand{\svskip}{\vspace{3mm}}
\newcommand{\C}{{\mathbb C}}
\newcommand{\BP}{{\mathbb P}}
\newcommand{\Q}{{\mathbb Q}}
\newcommand{\Z}{{\mathbb Z}}
\newcommand{\A}{{\mathbb A}}

\newcommand{\gc}{\EuFrak{c}}
\newcommand{\gm}{\EuFrak{m}}

\newcommand{\go}{\EuFrak{o}}
\newcommand{\gp}{\EuFrak{p}}

\newcommand{\SD}{{\mathcal D}}
\newcommand{\SK}{{\mathcal K}}
\newcommand{\SO}{{\mathcal O}}
\newcommand{\ST}{{\mathcal T}}

\newcommand{\Proof}{\noindent{\bf Proof.}\quad}
\newcommand{\Spec}{{\rm Spec}\:}
\newcommand{\Ker}{{\rm Ker}\:}
\newcommand{\Hom}{{\rm Hom}\:}
\newcommand{\SHom}{{\mathcal Hom}}
\newcommand{\Der}{{\rm Der}}
\newcommand{\SDer}{{\mathcal Der}}
\newcommand{\im}{{\rm Im}\:}

\newcommand{\height}{{\rm ht}\:}

\newcommand{\wt}{\widetilde}
\newcommand{\ol}{\overline}
\newcommand{\wh}{\widehat}
\newcommand{\st}[1]{\stackrel{{#1}}{\longrightarrow}}

\newcommand{\QED}{{\unskip\nobreak\hfil\penalty50\quad\null\nobreak\hfil
{$\Box$}\parfillskip0pt\finalhyphendemerits0\par\medskip}}

\newtheorem{thm}{Theorem}[section]
\newtheorem{lem}[thm]{Lemma}
\newtheorem{prop}[thm]{Proposition}
\newtheorem{cor}[thm]{Corollary}
\newtheorem{example}[thm]{Example}
\newtheorem{question}[thm]{Question}
\newtheorem{remark}[thm]{Remark}

\begin{document}
\title[Surjective derivations]{Surjective derivations in small dimensions}
\author{R. V. Gurjar, K. Masuda \and M. Miyanishi}
\thanks{The second and third authors are respectively supported by Grant-in-Aid for Scientific Research (C) 
15540043 and 21540055, JSPS} 
\address{School of Mathematics, Tata Institute of Fundamental Research \\
Homi Bhabha Road, Mumbai 400 001, INDIA}
\email{gurjar@math.tifr.res.in}
\smallskip

\address{School of Science and Technology \\
Kwansei Gakuin University \\
2-1 Gakuen, Sanda, 669-1337, JAPAN}
\email{kayo@kwansei.ac.jp}
\smallskip

\address{Research Center for Mathematical Sciences \\
Kwansei Gakuin University \\
2-1 Gakuen, Sanda, 669-1337, JAPAN}
\email{miyanisi@kwansei.ac.jp}
\date{July 24, 2012}

\maketitle

\begin{center}
\begin{small}{Dedicated to Professor C.S. Seshadri on his 80 th birthday}
\end{small}
\end{center}

\begin{abstract}
Let $D$ be a $\C$-derivation on a polynomial ring $\C[x_1,x_2]$. Cerveau \cite{cerveau} asserts that $D$ is surjective 
as a linear mapping if and only if $D=\frac{\partial}{\partial x_1}+ax_2\frac{\partial}{\partial x_2}$ with respect to 
a suitable algebraic change of coordinates of $\C[x_1,x_2]$, where $a \in \C$. Inspired by various results in \cite{cerveau}, 
we consider a surjective derivation defined on an affine domain over $\C$ of dimension one or two. Though our proofs 
are mostly algebraic or algebro-geometric, the idea using a result of Dimca-Saito \cite{DS} which is behind the arguments in 
\cite{cerveau} and based on the differential complex of the polynomial ring $\C[x_1,\ldots,x_n]$ is inspiring and affects 
our arguments.
\end{abstract}

\section*{Introduction} 
Let $X=\Spec B$ be an affine variety defined over the complex field $\C$ and let $\Omega^1_{X/\C}$ be the sheaf of 
Ka{\"h}ler differential $1$-forms on $X$. The tangent sheaf $\ST_{X/\C}$ is the dual sheaf of $\Omega^1_{X/\C}$, i.e., 
$\ST_{X/\C}=\SHom_{\SO_X}(\Omega^1_{X/\C},\SO_X)$ and it is identified with the sheaf of local derivations 
$\SDer_\C(\SO_X,\SO_X)$. A section $\theta \in \Gamma(X,\ST_{X/\C})$ is called a {\em vector field} on $X$. If $X$ is 
smooth and has dimension two, $\theta$ corresponds to an element of $\Gamma(X, \Omega^1_{X/\C}\otimes\SK^{-1}_{X/\C})$, 
where $\SK_{X/\C}=\Omega^2_{X/\C}$. In fact, the wedge product $\wedge : \Omega^1_{X/\C}\times\Omega^1_{X/\C} \to \Omega^2_{X/\C}$ 
gives rise to an isomorphism of $\SO_X$-modules $\Omega^1_{X/\C}\otimes \SK^{-1}_{X/\C} \cong \SHom_{\SO_X}(\Omega^1_{X/\C},\SO_X)$.  
The section $\theta$ corresponds to an element $\delta \in \Hom_B(\Omega^1_{B/\C}, B)$ and then to a $\C$-derivation 
$D \in \Der_\C(B,B)$ via $D=\delta\cdot d$, where $d : B \to \Omega^1_{B/\C}$ is the canonical differentiation. 

Let $P$ be a closed point of a smooth affine surface $X$ and let $\{x_1,x_2\}$ be a local system of parameters at $P$.
Then the derivation $D$ is expressed as 
$$
D=f_1\frac{\partial}{\partial x_1}+f_2\frac{\partial}{\partial x_2} \eqno{(\ast)}
$$
with $f_1, f_2 \in \SO_P$. We can associate a differential $1$-form $\omega \in \Omega^1_{X/\C,P}$ to $D$ by setting
$\omega=-f_2dx_1+f_1dx_2$. Then it is easy to see that 
\[
\omega\wedge dg=-D(g)dx_1\wedge dx_2
\]
for every $g \in \SO_P$. By a change of local coordinates $\{x'_1,x'_2\}$, we have 
\[
D=f'_1\frac{\partial}{\partial x'_1}+f'_2\frac{\partial}{\partial x'_2}, \quad \omega'=-f'_2dx'_1+f'_1dx'_2
\]
and $\omega'=J(x'/x)\omega$, where $J(x'/x)$ is the Jacobian determinant of $\{x'_1,x'_2\}$ with respect to 
$\{x_1,x_2\}$. Hence the correspondence $g \mapsto D(g)$ is independent of the choice of local coordinates. 
The derivation $D$ is {\em surjective} if it is surjective as an endomorphism of the set $B$. If $B=\C[x_1,x_2]$ is 
a polynomial ring, the derivation $D$ of $B$ is surjective if and only if, for every $f \in B$, there exists an 
element $g \in B$ such that $\omega\wedge dg=-fdx_1\wedge dx_2$, where $\omega$ is defined on $B$ by $(\ast)$.

In a very interesting article \cite{cerveau}, Cerveau asserts as stated in the abstract that a $\C$-derivation $D$ of 
a polynomial ring $\C[x_1,x_2]$ is surjective if and only if 
\[
D=\frac{\partial}{\partial x_1}+ax_2\frac{\partial}{\partial x_2}  
\]
after an {\em algebraic} change of variables, where $a \in \C$. Unfortunately, there exists an error in the proof 
of \cite[Lemma 5.5]{cerveau} which is not corrected for the moment and thereby suspends the validity of the 
above-mentioned result. Furthermore, an analytic theory, e.g., the properties of analytic foliations associated 
to a differential $1$-form is used occasionally in the proof of \cite{cerveau}. 

Inspired by Cerveau's work and his idea, we intend to generalize the results to surjective derivations $D$ of affine 
domains $B$ over $\C$ with small dimension. Our trial is successful only in the special cases. An element $f$ of $B$ 
is called an {\em integral} element with respect to $D$ if $D(f)$ is divisible by $f$ in $B$. In Nowicki \cite{A2}, 
an integral element for a derivation on a polynomial ring is called a {\em Darboux polynomial}. 

In the present article, we prove the following two results Theorems $1$ and $2$.
\svskip

\noindent
{\bf Theorem 1.}\ \ {\em 
Let $D$ be a $\C$-derivation of $\C[x_1,x_2]$. Then $D={\partial}/{\partial x_1}$ with respect to a set of variables 
$\{x_1,x_2\}$ if and only if $D$ is surjective and $\Ker D \supsetneqq \C$.}
\svskip

We have only to prove the {\em if} part because the {\em only if} part is obvious. Since a locally nilpotent 
$\C$-derivation on $\C[x_1,x_2]$ is always written as $f(x_2){\partial}/{\partial x_1}$ with respect to a suitable set 
of variables $\{x_1,x_2\}$, the above result gives a characterization for a given reduced $\C$-derivation on 
$\C[x_1,x_2]$ to be locally nilpotent, where a derivation is {\em reduced} if the coefficients $f_1, f_2 \in \C[x_1,x_2]$ 
have no irreducible common factors when we write 
\[
D=f_1\frac{\partial}{\partial x_1}+f_2\frac{\partial}{\partial x_2}.
\] 

We prove the above theorem in a generalized setting where $\C[x_1,x_2]$ is replaced by a factorial affine domain $B$
of dimension two over $\C$ such that $B^*=\C^*$. In this case, the reducedness of $D$ is equivalent to the 
condition that the divisorial part of $D$ on the smooth part of $\Spec B$ is zero (cf. \cite {MM1}). Thus we 
obtain a different algebraic characterization of $\C[x_1,x_2]$ in terms of a surjective derivation. 

In the subsequent proof, we do not make a full use of the surjectivity of $D$ but only a partial condition 
derived from the surjectivity of $D$, that is, the existence of an element $t \in B$ with $D(t)=\delta(dt)=1$. 
If $R$ is a normal affine domain of dimension one, we take $R$ as the above $B$. Then the module of differential forms 
$\Omega^1_{R/\C}$ is a free $R$-module $Rdt$ since the $R$-module homomorphism $\delta : \Omega^1_{R/\C} \to R$ induces 
an isomorphism, and $\delta(\omega)=g \in R$ for every $\omega \in \Omega_{R/\C}$ if we write $\omega=gdt$. 
Since $D$ is surjective, $g=D(h)$ for some element $h \in R$. Hence it follows that $\delta(\omega)=D(h)=\delta(dh)$, 
i.e., $\delta(\omega-dh)=0$. Since $\delta$ is an isomorphism of $R$-modules, we have $\omega=dh$. Namely, every (closed) 
form on a curve $\Spec R$ is {\em exact} in the case of $\dim R=1$. In view of the isomorphism $\delta$ for $R$, 
the exactness of any differential $1$-form is equivalent to the surjectivity of $D$. In the case $\dim B \ge 2$, 
see \cite[Proposition 1.2]{cerveau} for the interpretation of $D$ being surjective. We then use the mixed Hodge structure 
on $\Spec R$ to conclude that $\Spec R$ is in fact the affine line. In the case where $B=\C[x_1,x_2]$, the assumption that 
$\Ker D \supsetneqq \C$ yields an $\A^1$-fibration $\Spec B \to \Spec A$ with $A=\Ker D$ and consequently gives a proof 
of the above theorem 1. One can say that the surjectivity of $D$ is {\em less} algebraic than it appears from 
the definition. We give several proofs of Theorem \ref{Theorem 1.2} when $\dim R=1$ in order to illustrate the condition 
of surjectivity of $D$. It is still a mysterious condition when $\dim B > 1$. 

Although there is a gap in the proof, the correct parts combined together in \cite{cerveau} imply the following result.
\svskip

\noindent
{\bf Theorem 2.}\ \ {\em Let $D$ be a surjective derivation of $\C[x_1,x_2]$. Then there exists an integral element with 
respect to $D$.}
\svskip

It should be remarked that most $\C$-derivations of $\C[x_1,x_2]$ have no integral elements. In Nowicki \cite{A, A2}, 
a $\C$-derivation without integral elements is called {\em simple} and many examples are given. 

Besides the reference \cite{cerveau}, some observations have been made in \cite{Essen, Stein} on surjective derivations 
on a polynomial ring $R[x,y]$ of dimension two over $R=\C$ or $R$ being a $\Q$-algebra, especially on the relation 
between surjectivity and local nilpotence of $D$ under the assumption that $D$ has divergence zero. In fact, if $D$ is 
a surjective derivation of $\C[x,y]$ with divergence zero (see \cite{Essen,F} for the definition) then $\Ker D \supsetneqq 
\C$. Hence the above theorem implies that $D$ is locally nilpotent. This is proved in \cite{Essen} for a $\Q$-algebra 
$R$ instead of $\C$. 

When we looked for examples of derivations without integral elements, Nina Gupta of TIFR kindly found us the references of Nowicki
\cite{A,A2}. L. Makar-Lomanov also constructed for us an example (see Example \ref{Example 3.13}). We are very 
thankfulful to 
both of them.

\section{Case of dimension one and related remarks}

Let $R$ be a normal affine domain of dimension one. Suppose that $R$ has a surjective $\C$-derivation. Let $t$ 
be an element of $R$ such that $D(t)=1$. Let $R_0=\C[t]$. Then $R_0$ is a subring of $R$ such that 
$D(R_0) \subseteq R_0$ and $D$ is locally nilpotent on $R_0$. The derivation $D$ corresponds to an element 
$\delta \in \Hom_R(\Omega_{R/\C},R)$ such that $D(z)=\delta(dz)$ for every $z \in R$, where $\delta$ is 
an $R$-homomorphism. Since $D(t)=1$, $\delta$ is a surjective homomorphism. Consider an exact sequence of 
differential $R$-modules
\[
\Omega_{R_0/\C}\otimes_{R_0}R \st{\alpha} \Omega_{R/\C} \st{\beta} \Omega_{R/R_0} \to 0,
\]
where $\Omega_{R_0/\C}\otimes_{R_0}R$ is a free $R$-module of rank one generated by $dt\otimes 1$ and hence the 
homomorphism $\alpha$ is injective. Since $\Spec R$ is smooth by the hypothesis, $\Omega_{R/\C}$ is a locally 
free module of rank one and the homomorphism $\delta$ is an $R$-isomorphism since it is surjective. In particular, 
$\Omega_{R/\C}$ is an $R$-free module generated by $dt$. Then $\Omega_{R/R_0}=0$, which implies that the induced 
morphism $q : \Spec R \to \Spec R_0$ is \'etale. 

\begin{lem}\label{Lemma 1.1}
The morphism $q$ is surjective. Hence a linear form $t-c$ is not invertible in $R$ for every $c \in \C$.
\end{lem}
\Proof
Suppose that $q$ is not surjective and the point $t=0$ is not in the image of $q$. Then $t$ is invertible in $R$.
Hence there exists an element $u$ of $R$ such that $D(u)=1/t$. Write the minimal equation of $u$ over $\C(t)$ in 
the form 
\begin{eqnarray*}
&& a_0(t)u^n+a_1(t)u^{n-1}+\cdots+a_{n-1}(t)u+a_n(t)=0,\\
&& \quad \forall~ a_i(t) \in R_0,~\gcd(a_0(t),\ldots,a_n(t))=1.  
\end{eqnarray*}
If we have another equation of $u$ with degree $n$
\[
b_0(t)u^n+b_1(t)u^{n-1}+\cdots+b_{n-1}(t)u+b_n(t)=0,~~\forall~ b_i(t) \in R_0,
\]
the condition $\gcd(a_0(t),\ldots,a_n(t))=1$ implies $b_i(t)/a_i(t)$ is an element of $R_0$ which is independent of 
$i$. Applying $D$ to the first equation, we have 
\begin{eqnarray*}
&& ta'_0(t)u^n+\left(ta'_1(t)+na_0(t)\right)u^{n-1}+\cdots \\
&&+\left(ta'_i(t)+(n-i+1)a_{i-1}(t)\right)u^{n-i}+\cdots+\left(ta'_n(t)+a_{n-1}(t)\right)=0.
\end{eqnarray*}
Hence we have
\[
\frac{ta'_0(t)}{a_0(t)}=\frac{ta'_1(t)+na_0(t)}{a_1(t)}=c \in \C. 
\]
Let $r=\deg a_0(t)$. Then the equality $ta'_0(t)=ca_0(t)$ implies that $c=r$ and $a_0(t)=c_0t^r$ with $c_0 \in \C^*$. 
We may assume that $c_0=1$. Write
\[
a_1(t)=d_0t^s+d_1t^{s-1}+\cdots+d_{s-1}t+d_s, \quad d_0 \ne 0.
\]
Then the equation $ta'_1(t)+na_0(t)=ca_1(t)$ is written as
\[
sd_0t^s+(s-1)d_1t^{s-1}+\cdots+d_{s-1}t+nt^r=rd_0t^s+rd_1t^{s-1}+\cdots+rd_{s-1}t+rd_s.
\]
This implies that $s=r$ and hence $n=0$. This is a contradiction. 
\QED

The following result determines the structure of $R$.

\begin{thm}\label{Theorem 1.2}
Let $R$ be a normal affine domain of dimension $1$ over $\C$. Assume that $R$ has a surjective $\C$-derivation $D$. Then the above 
morphism $q$ is an isomorphism. Hence $R=\C[t]$ and $D=d/dt$.
\end{thm}
\Proof
Let $X=\Spec R$, let $C$ be a smooth completion of $X$ and let $\Delta=C\setminus X$, which we regard as a reduced effective 
divisor on $C$. The mixed Hodge structure theorem due to Deligne (see \cite[Chapter 8]{V}) gives the following isomorphism 
\[
H^1(X;\C) \cong H^0(C,K_C+\Delta)\oplus H^1(C, \SO_C),
\]
where we note that $\Omega^1_{C/\C}(\log \Delta)=\SO_C(K_C+\Delta)$. By the Riemann-Roch theorem, we have 
\[
h^0(C,K_C+\Delta)=g-1+n,~~h^1(C,\SO_C)=g,~~n=\deg \Delta.
\]
Hence $\dim H^1(X;\C)=2g-1+n$. The de Rham theorem gives a non-degenerate pairing 
$$
H_1(X;\C)\times H^1(X;\C) \to \C, \quad (\gamma, \omega) \mapsto \int_\gamma \omega, \eqno{(1)}
$$
where $\gamma$ is a closed loop in $X$ and $\omega$ is a differential $1$-form having at most logarithmic singularities 
along $\Delta$. Then $\omega$ restricted onto $X$ is a regular form, and hence $\omega|_X=dh$ for an element $h \in R$. 
This follows from the surjectivity of $D$. Since $\gamma$ is a loop in $X$, we have 
\[
\int_\gamma \omega=\int_\gamma \omega|_X=\int_\gamma dh=0.
\]
The last equality follows from Stokes' theorem. Since the pairing (1) is non-degenerate, this implies that $H^1(X;\C)=0$.
Namely, $2g-1+n=0$. Since $n > 0$, this is equivalent to $g=0$ and $n=1$. Hence $X$ is isomorphic to the affine line. 

Now write $R=\C[z]$. Then $D=d/dz$ possibly after a suitable change of $z$. Since $D(t)=1$, we have $t=z+c$ with $c \in \C$. 
Hence $R=\C[t]$. 
\QED

\noindent
{\em The second proof of Theorem \ref{Theorem 1.2}.}\ \ In the above proof, one can use the algebraic de Rham theorem 
by Grothendieck \cite{G}. In fact, if $X$ is a smooth affine variety, it is shown that the following equality holds
\[
H^*(X;\C)=H^*(X,\Gamma(X,\Omega^*_{X/\C})).
\]
Hence if $\dim X=1$ and every regular $1$-form on $X$ is exact, it follows that $H^1(X;\C)=0$, which implies that 
$X \cong \A^1$.
\QED
\svskip

\noindent
{\em The third proof of Theorem \ref{Theorem 1.2}.}\ \ We can give a different proof without using the de Rham theorem. 
In the proof of Lemma \ref{Lemma 1.1}, it is shown that the morphism $q : \Spec R \to \Spec R_0$ is surjective and 
\'etale. Let $X=\Spec R$ and $Z=\Spec R_0$. Since $Z$ is simply connected as $Z \cong \A^1$, $q$ is an isomorphism if 
we show that $q$ is finite. This follows from the following lemma. In fact, let $\wt{R}$ be the normalization of $R_0$ 
in $R$ and let $\wt{X}=\Spec \wt{R}$. Then $X$ is an open set of $\wt{X}$. If $X \subsetneqq \wt{X}$, then $R$ cannot have a 
surjective derivation. Hence $X=\wt{X}$ and $q$ is finite. 

\begin{lem}\label{Lemma 1.3}
Let $C$ be a smooth affine curve and let $P$ be a point of $C$. Then there exists a regular $1$-form on $C$ which is not exact 
on $C\setminus\{P\}$.
\end{lem}
\Proof
Let $\Omega^1_{C/\C}(P)$ be the sheaf of $1$-forms with at most a simple pole at $P$. Then $\Omega^1_{C/\C}(P)$ is a 
coherent $\SO_C$-Module containing properly $\Omega^1_{C/\C}$ as an $\SO_C$-Module. Since $C$ is affine, we have a 
proper containment
\[
\Gamma(C,\Omega^1_{C/\C}(P)) \supsetneqq \Gamma(C,\Omega^1_{C/\C}).
\]
Hence there exists a $1$-form $\eta$ on $C$ which is regular on $C\setminus\{P\}$ and has a simple pole at $P$. 
Suppose that $\eta$ is exact, i.e., $\eta=dh$ for a rational function $h$. Then $h$ must be regular. In fact, if $h$ 
has a pole of order $n$ at $P$, then $dh$ has a pole of order $n+1$. Hence $\eta$ cannot have a simple pole at $P$. 
But this contradicts the choice of $\eta$. 
\QED

Employing the argument in the proof of Lemma \ref{Lemma 1.3}, we can show the following result.

\begin{lem}\label{Lemma 1.4}
Let $f:X \to Z$ be a non-constant morphism between smooth affine curves. Let $\eta$ be a regular $1$-form on $Z$.  
If $f^*(\eta)$ is exact on $X$, then $\eta$ is exact on $Z$.
\end{lem}
\Proof 
Let $\wt{f} : \wt{X} \to Z$ be the normalization of $Z$ in the function field of $X$. Then $X$ is an open subset of $\wt{X}$. 
By the assumption, $f^*(\eta)=dh$ for a regular function $h$ on $X$. The proof of Lemma \ref{Lemma 1.3} shows that 
it is exact on $\wt{X}$. In fact, if $h$ has a pole of order $n$ at a point $P$ of $\wt{X}\setminus X$, then $dh$ has 
a pole of order $n+1$ at $P$ which contradicts the regularity of $f^*(\eta)$. 
Let $W$ be the normalization of $Z$ in the smallest Galois extension of the function field of $Z$ containing 
the function field of $\wt{X}$. Then the normalization morphism $\nu : W \to Z$ splits as $\nu=\wt{f}\circ \mu$, where 
$\mu : W \to \wt{X}$ is also a Galois extension. Let $G$ be the Galois group of the extension $\C(W)/\C(Z)$. 
Then $\nu^*(\eta)=\mu^*(\wt{f}^*(\eta))=\mu^*(dh)=d(\mu^*(h))$. Let $g \in G$ and let $\alpha(g) : W \to W$ be the 
automorphism induced by the field automorphism $g : \C(W) \to \C(W)$. Then we have
\[
\alpha(g)^*(\nu^*(\eta))=d(g(\mu^*(h))),
\]
where $\alpha(g)^*(\nu^*(\eta))=\nu^*(\eta)$ for every $g \in G$. Hence by the averaging trick, we have
\[
\eta=\frac{1}{|G|}\sum_{g\in G}d(g(\mu^*(h)))=d\bigg(\frac{1}{|G|}\sum_{g\in G}g(\mu^*(h))\bigg), 
\]
where $\frac{1}{|G|}\sum_{g\in G}g(\mu^*(h))$ is a regular function on $Z$. So, $\eta$ is exact on $Z$.
\QED

Lemma \ref{Lemma 1.4} gives a different proof of Lemma \ref{Lemma 1.1}.
\svskip

\noindent
{\em The second proof of Lemma \ref{Lemma 1.1}}\ \ With the notations in Lemma \ref{Lemma 1.1}, suppose that 
$q : X \to \A^1$ is not surjective. By Lemma \ref{Lemma 1.3}, we can find a regular $1$-form $\eta$ on $q(X)$ 
which is not exact on $q(X)$. By Lemma \ref{Lemma 1.4}, $q^*(\eta)$ is not exact on $X$. This is a contradiction 
since every regular $1$-form on $X$ is exact by the surjectivity of $D$.
\QED

To understand what the surjectivity of $D$ implies, we conduct an elementary calculation to verify Theorem 
\ref{Theorem 1.2} in the case where $R$ is rational over $\C$. 

\begin{lem}\label{Lemma 1.5}
With the above notations, suppose that $R$ is rational and $D$ is surjective on $R$. Then $R=\C[t]$ and $D=d/dt$.
\end{lem}
\Proof
We can write
\[
R=\C\left[t, \frac{1}{\prod_{i=1}^r(t-\alpha_i)}\right],\quad \alpha_i \in \C,~~\alpha_i \ne \alpha_j~(i \ne j).
\]
We only consider the case where the derivation $D$ is written as 
\[
D=\frac{f(t)}{\prod_{i=1}^r(t-\alpha_i)^{n_i}}\frac{d}{dt}, \quad f(t) \in \C[t],~~f(\alpha_i) \ne 0~(\forall~i),~~ n_i > 0.
\]
Since $D$ is surjective, there exists an element $u \in R$ such that $D(u)=1$. Write
\[
u=\frac{g(t)}{\prod_{i=1}^r(t-\alpha_i)^{m_i}}, \quad g(\alpha_i) \ne 0~(\forall~i),~~m_i \in \Z.
\]
The equation $D(u)=1$ yields
\[
\frac{f(t)}{\prod_{i=1}^r(t-\alpha_i)^{n_i}}\cdot\left\{\frac{g'(t)}{\prod_{i=1}^r(t-\alpha_i)^{m_i}} 
-\frac{g(t)}{\prod_{i=1}^r(t-\alpha_i)^{m_i}}\cdot \sum_{i=1}^r\frac{m_i}{(t-\alpha_i)}\right\}=1 .
\]
This is written as follows
\begin{eqnarray*}
&& f(t)\left\{g'(t)\prod_{i=1}^r(t-\alpha_i)- g(t)\sum_{i=1}^r m_i\cdot\frac{\prod_{i=1}^r(t-\alpha_i)}{t-\alpha_i} \right\}  \\
&&\qquad = \prod_{i=1}^r(t-\alpha_i)^{n_i+m_i+1}.
\end{eqnarray*}
Since $f(\alpha_i) \ne 0$ for all $i$, it follows that $n_i+m_i+1=0$ for all $i$ and hence $f(t) \in \C^*$. After 
changing $g(t)$ by $f(0)g(t)$, we may assume that $f(t)=1$. Then $\deg g(t) =s > 0$, and $g(t)=0$ has no multiple roots. 
We can write 
\[
g(t)=c\prod_{j=1}^s(t-\beta_j), ~~\beta_j \ne \beta_\ell~(j \ne \ell),~~\beta_j \ne \alpha_i.
\]
The above equation is then written as 
\[
\sum_{j=1}^s\frac{c}{t-\beta_j}-\sum_{i=1}^r\frac{cm_i}{t-\alpha_i}=\frac{1}{\prod_{i=1}^r(t-\alpha_i)\prod_{j=1}^s(t-\beta_j)},
\]
where $-m_i=n_i+1 \ge 2$. Assume that $r \ge 1$. Then comparison of the top $t$-degrees of both sides of this equation 
after multiplied $\prod_{i=1}^r(t-\alpha_i)\prod_{j=1}^s(t-\beta_j)$ yields an equality
\[
1-\sum_{i=1}^rm_i=0.
\]
But this is a contradiction. If $r=0$, then $R=\C[t]$. 
\QED

The use of Zariski's lemma \cite[Lemma 4]{Z2} enables us to drop in the statement of Theorem \ref{Theorem 1.2} the assumption 
that $R$ is normal. 

\begin{thm}\label{Theorem 1.6}
Let $(R,\gm)$ be a geometric local ring \footnote{Namely, $R$ is the local ring of a closed point of an algebraic 
variety.} and let $D$ be a $\C$-derivation of $R$. Assume that there exists an element $t \in R$ with $D(t)=1$. 
Then the following assertions hold.
\begin{enumerate}
\item[(1)]
If $\dim R=1$ then $R$ is normal.
\item[(2)]
Assume that $\dim R=2$ and $R$ is normal. Then $R$ is regular.
\end{enumerate}
\end{thm}
\Proof
Let $\wh{R}$ be the $\gm$-adic completion of $R$. Then the derivation $D$ extends to a derivation on $\wh{R}$ 
in a natural fashion, which we denote by $\wh{D}$. Then we have $\wh{D}(t)=1$. By Zariski's lemma,  
$\wh{R}=\go[[\:t\:]]$ for a complete local ring $\go$ such that the restriction of $\wh{D}$ onto $\go$ is zero and 
$t$ is analytically independent over $\go$. Note that $\wh{R}$ is reduced by Nagata \cite[Theorem (37.5)]{N}. 
If $\dim R=1$ then $\dim \go=0$ and hence $\go=\C$. So, $\wh{R}=\C[[\:t\:]]$, whence $R$ is normal because $\wh{R}$ 
is a faithfully flat $R$-module. If $\dim R=2$ and $R$ is normal, then $\wh{R}$ is normal by \cite[Theorem (37.5)]{N}. 
Hence $\go$ is a complete normal ring of dimension one. Then $\go=\C[[\:u\:]]$ and $\wh{R}=\C[[t,u]]$, whence $R$ is 
regular by the faithfully flat descent of regularity. 
\QED

Note that the hypothesis that $D$ is surjective is stronger than the hypothesis that there exists an element $t \in R$ 
with $D(t)=1$. If we use the stronger hypothesis, we can give a different proof for $R$ being normal if $\dim R=1$. 

\begin{prop}\label{Proposition 1.7}
Let $R$ be an affine domain of dimension one defined over $\C$. Assume that $R$ has a surjective $\C$-derivation $D$. 
Then $R$ is normal.
\end{prop}
\Proof
Suppose that $R$ is not normal. Let $\wt{R}$ be the normalization of $R$. Then, by \cite{S}, the derivation $D$ extends to 
a $\C$-derivation $\wt{D}$ of $\wt{R}$. Let $\gc$ be the conductor of $R$ in $\wt{R}$. Then $\gc$ is an ideal of $R$ and 
$\wt{R}$ as well. Furthermore, the closed set $V(\gc)$ is the singular locus of $\Spec R$ and hence $R/\gc$ is 
a $\C$-vector space of finite dimension. The conductor $\gc$ is defined as the set 
\[
\gc=\{b \in \wt{R} \mid bx \in R,~\forall~x \in \wt{R}\}. 
\]
Let $b \in \gc$. Then, for every $x \in \wt{R}$, we have $D(bx)=D(b)x+b\wt{D}(x)$, whence $D(b)x \in R$ as 
$D(bx), b\wt{D}(x) \in R$. So, $D(\gc) \subseteq \gc$. Hence $D$ induces a $\C$-linear endomorphism $\ol{D}$ of $R/\gc$. 
We have the following commutative exact diagram of $\C$-modules
\[
\begin{array}{lllllllll}
0 & \to & \gc & \to & R & \to & R/\gc & \to & 0 \\
\mbox{}\\
 & & \Big\downarrow D_c & & \Big\downarrow D & & ~\:\Big\downarrow \ol{D} & & \\
 \mbox{}\\
0 & \to & \gc & \to & R & \to & R/\gc & \to & 0 
\end{array}
\]
where the middle vertical arrow $D$ is a surjection and the right vertical arrow $\ol{D}$ is thereby a surjective 
$\C$-linear endomorphism of the finite-dimensional $\C$-vector space $R/\gc$. Hence $\ol{D}$ is a bijection. 
Then, by the snake lemma, $\Ker D_c=\Ker D$, and $\Ker D = \C$ for $D$ is a non-zero $\C$-derivation of $R$ and $\dim R=1$. 
This is a contradiction because $\gc\cap \C=(0)$. This shows that $R$ is normal. 
\QED

The surjectivity of $D$ does not necessarily imply the normality in the higher-dimensional case.

\begin{example}\label{Example 1.8}{\em
Let $B=\C[X,Y,Z]/(X^2+Y^3)$ and let $x,y,z$ be respectively the residue classes of $X,Y,Z$ in $B$. Define a derivation
$\SD$ on $\C[X,Y,Z]$ by 
\[
\SD=\frac{1}{2}X\frac{\partial}{\partial X}+\frac{1}{3}Y\frac{\partial}{\partial Y}+\frac{\partial}{\partial Z}.
\]
Then $\SD$ maps the ideal $(X^2+Y^3)$ into itself. Hence $\SD$ induces a $\C$-derivation $D$ on the ring $B$. Since 
$\SD$ is easily verified to be surjective (cf. Lemma \ref{Lemma 1.11} below), the derivation $D$ is surjective as well. 
The ring $B$ is not normal.}
\end{example}

In the case of an affine plane curve, we have the following result.

\begin{prop}\label{Proposition 1.9}
Let $R:=\C[X,Y]/(F)$ be the coordinate ring of an irreducible, affine plane curve $C:=\{F(X,Y)=0\}$. Let $x, y$ be 
the residue classes of $X, Y$ in $R$ respectively. Write the residue classes of the partial derivatives $F_X, F_Y$ 
by $f_x, f_y$ respectively. Then the following assertions hold.
\begin{enumerate}
\item[(1)]
There exists a $\C$-derivation $D$ of $R$ with $D(x)=1$ if and only if $f_x \in f_yR$. If this condition is 
satisfied, then the curve $C$ is smooth.
\item[(2)]
Suppose that there exists a $\C$-derivation $D$ of $R$ with $D(x)=1$. Then either $R$ is a polynomial ring $\C[x]$ 
or the element $f_y$ is a non-constant, invertible element in $R$.
\end{enumerate}
\end{prop}
\Proof
(1)\ Let $\Delta$ be a $\C$-derivation on $\C[X,Y]$. Then $\Delta=A(\partial/\partial X)+B(\partial/\partial Y)$ 
with $A, B \in \C[X,Y]$, and $\Delta$ induces a $\C$-derivation $D$ on $R$ if and only if $\Delta(F) \in (F)$. 
Further, if $D$ is the induced derivation, then we have $f_xD(x)+f_yD(y)=0$. If $D(x)=1$, then $f_x =-f_yD(y) 
\in f_yR$ as $D(y) \in R$. Conversely, suppose that $f_x=-f_yg$ with $g \in R$. Then there exist elements 
$G, H \in \C[X,Y]$ such that $g$ is the residue class of $G$ and $F_X=-F_YG+FH$. Set $\Delta=(\partial/\partial X)
+G(\partial/\partial Y)$. Then $\Delta(F)=F_X+F_YG=FH$. Hence $\Delta$ induces a $\C$-derivation $D$ with $D(x)=1$.
By Theorem \ref{Theorem 1.6}, $R$ is normal and hence the curve $C$ is smooth if this condition is satisfied. 

(2)\ Since $C$ is smooth, the ideal $(f_x,f_y)$ is a unit ideal. Suppose that $f_y$ is a constant $c \in \C$. 
Then $F_Y-c$ is divisible by $F$. Since $\deg_YF_Y < \deg_YF$, it follows that $F_Y=c$. This implies that $F=cY-G(x)$.
So, $R=\C[x]$. Suppose that $f_y$ is not a constant. Since $(f_x,f_y)=R$, we find elements $a,b \in R$ such that 
$af_x+bf_y=1$ in $R$. Hence $af_x/f_y+b=1/f_y$. By the assertion (1), $f_x/f_y \in R$ and hence $1/f_y \in R$.
\QED

By the following remark, there exist an affine normal domain $R$ over $\C$ with genus $g > 0$ and a $\C$-derivation 
$D$ such that $D(x)=1$ for some element $x \in R$. 

\begin{remark}\label{Remark 1.10}{\em 
Let $C$ be a smooth projective curve of genus $g$. By a classical result of F. Severi (cf. \cite[Proposition 8.1]{Fulton}), 
there is a morphism $f : C \to \BP^1$ of degree $n \ge g+1$ such that, over every point $P \in \BP^1$, there is at most 
one point, say $\wt{P}$ (if it exists) which is ramified over $P$ and has ramification index $2$. Fix a point 
$P_{\infty} \in \BP^1$. Let $C_0$ be the open set obtained from $C$ by omitting all the points lying over $P_{\infty}$ 
and also all the ramified points. Let $R$ be the coordinate ring of the affine curve $C_0$. The induced morphism 
$f : C_0 \to \A^1$ is surjective and \'etale. Let $\C[x]$ be the coordinate ring of $\A^1$. Since $dx$ is a 
nowhere vanishing regular $1$-form on $C_0$, we have $\Omega_{R/\C}=Rdx$. In fact, for every element $a \in R$, the 
ratio $da/dx$, which lies in the function field of $R$, has non-negative value at the discrete valuation ring 
corresponding to any point of $C_0$, i.e., $da/dx \in R$. This implies that $D:=d/dx$ is a $\C$-derivation of $R$ 
and that $D(x)=1$. We have seen in Theorem \ref{Theorem 1.2} that if $g>0$ then $D$ cannot be surjective.
Furthermore, since $\Omega^1_{R/\C}=Rdx$, $R$ is a complete intersection by a result of Murthy-Towber \cite{MT,K}. 
Hence there are many examples of normal $1$-dimensional affine domains with a $\C$-derivation $D$ such that 
$D(x)=1$ for some $x \in R$ but $D$ is not surjective.} 
\QED
\end{remark}

The following result gives a construction of a surjective derivation. A $\C$-derivation $D$ on a $\C$-algebra $R$ is 
called {\em locally finite} if, for every element $a \in R$, the $\C$-vector space spanned by the set $\{D^i(a)\mid 
i \ge 0\}$ has finite dimension.

\begin{lem}\label{Lemma 1.11}
Let $R$ be a $\C$-algebra and let $D$ be a $\C$-derivation on $R$. Let $t$ be an indeterminate and let $B=R[t]$. 
Define a $\C$-derivation $\wt{D}$ on $B$ by $\wt{D}(a)=D(a)$ for every $a \in R$ and $\wt{D}(t)=1$. Then $\wt{D}$ 
is surjective if $D$ satisfies one of the conditions:
\begin{enumerate}
\item[(1)]
$D=0$.
\item[(2)]
$D$ is locally nilpotent.
\item[(3)]
$D$ is locally finite.
\item[(4)]
$R\setminus \C \subseteq \im D$. 
\end{enumerate}
\end{lem}
\Proof
Since we have implications $(2) \Rightarrow (1)$ and $(2) \Rightarrow (3)$, we prove the assertion in the cases (3) and (4). 

(3)\ Let $a$ be an element of $R$. It suffices to show that $at^n \in \im \wt{D}$ for every $n \ge 0$. 
If $D(a)=0$ then $at^n=\wt{D}(at^{n+1})/(n+1) \in \im \wt{D}$. Suppose that $D^m(a)=0$. By induction on $m$, 
we show that $at^n \in \im \wt{D}$ for every $n \ge 0$. In fact, we have
\[
\wt{D}(at^{n+1})=(n+1)at^n+D(a)t^{n+1},
\]
where $D(a)t^{n+1} \in \im \wt{D}$ by the induction hypothesis because $D^{m-1}(D(a))$ $=0$. Hence $at^n \in \im \wt{D}$. 
Meanwhile, since $D$ is locally finite, there exists a relation
\[
c_0a+c_1D(a)+\cdots+c_pD^p(a)=0, \ \ \forall~c_i \in \C.
\]
If $c_0 \ne 0$, then $a \in \im D$. If $c_0=\cdots=c_{\ell-1}=0$ and $c_\ell\ne 0$, then $a-a_0\in \im D$ for an 
element $a_0 \in R$ with $D^\ell(a_0)=0$. In fact, it suffices to put 
\[
a_0=a+ \frac{1}{c_\ell}\big(c_{\ell+1}D(a)+\cdots+c_pD^{p-\ell}(a)\big).
\]
Since $a_0t^n \in \im \wt{D}$ for $\forall~n \ge 0$ by the above observation, we may assume $a \in \im D$ 
in order to show that $at^n \in \im \wt{D}$. Write $a=D(a')$ with $a' \in R$. Since 
\[
at^n=D(a')t^n=\wt{D}(a't^n)-na't^{n-1},
\]
we are done by induction on $n$.

(4)\ Note that $ct^n \in \im\wt{D}$ for $c \in \C$ and $n \ge 0$, for $\wt{D}(ct^{n+1}/(n+1))=ct^n$. If $a \in 
R\setminus\C$, we show by induction on $n$ that $at^n \in \im\wt{D}$. If $n=0$ then $a \in \im D$ by the 
assumption. Write $a=D(a')$ with $a' \in R$. Then $at^n=D(a')t^n=D(a't^n)-na't^{n-1}$ and $a't^{n-1} \in \im\wt{D}$ 
by the induction hypothesis, whence $at^n \in \im\wt{D}$.
\QED

\section{Proof of Theorem 1}

Let $B$ be a factorial affine domain of dimension two over $\C$ such that $B^*=\C^*$. Throughout the present section, 
the domain $B$ satisfies these conditions. We assume that $B$ is endowed with a $\C$-derivation $D$ such that $\Ker D 
\supsetneqq \C$. Wherever we need the surjectivity condition on $D$, we mention it explicitly. 

\begin{lem}\label{Lemma 2.1}
Let $D$ be a non-zero $\C$-derivation on $B$ such that $\Ker D \supsetneqq \C$. Let $A=\Ker D$. Then the 
following assertions hold.
\begin{enumerate}
\item[(1)]
$A$ is a polynomial ring in one variable.
\item[(2)]
$Q(A)$ is algebraically closed in $Q(B)$, where $Q(A)$ and $Q(B)$ are the quotient fields of $A$ and $B$ respectively.
\end{enumerate}
\end{lem}
\Proof
(1)\ The ring $A$ is integrally closed in $B$. In fact, if $z$ is an element of $B$ which is integral over $A$, let
\[
f(z)=z^n+a_1z^{n-1}+\cdots+a_{n-1}z+a_n=0, \quad \forall a_i \in A
\]
be a monic equation of $z$ over $A$. We may assume that the degree $n$ is minimal among such monic equations for $z$. 
Applying $D$ to $f(z)$, we obtain a relation 
\[
\frac{1}{n}f'(z)D(z)=0,
\]
where $\dfrac{1}{n}f'(z)$ is monic in $z$ and hence nonzero by the hypothesis. Then $D(z)=0$. So, $z \in A$. It is clear 
that $\dim A=1$ as $D$ is non-trivial and $A=B\cap \Ker \wt{D}$, where $\wt{D}$ is the $\C$-derivation of $Q(B)$ 
which is the natural extension of $D$. Hence $A$ is finitely generated by a lemma of Zariski \cite{Z}. Furthermore, 
$Q(A)$ is rational over $\C$ because $B$ is factorial. Since $A^*=\C^*$ as $B^*=\C^*$, it follows that $A$ is 
a polynomial ring $\C[\xi]$ in one variable $\xi$.

(2)\ Suppose that $z$ is an element of $Q(B)$ which is algebraic over $Q(A)$. Then there is an algebraic relation 
\[
a_0z^n+a_1z^{n-1}+\cdots+a_{n-1}z+a_n=0, ~~ \forall a_i \in A,~a_0 \ne 0.
\]
Then the element $a_0z$ is integral over $A$. Hence $a_0z \in A$ by the assertion (1). So, $z \in Q(A)$ 
and $Q(A)$ is algebraically closed in $Q(B)$.
\QED

Let $Y=\Spec B$, $X=\Spec A$ and $p : Y \to X$ the morphism associated to the natural inclusion 
$A \hookrightarrow B$. Since $Q(A)$ is algebraically closed in $Q(B)$, the general fibers of $p$ are smooth curves. 
The derivation $D$ induces a derivation $D_c$ on the fiber $F_c=\Spec B/(\xi-c)$, which is surjective provided
so is $D$.

We consider, in general, a fibration $p : Y \to X$, which is by definition a dominant morphism of algebraic varieties 
$Y, X$ with smooth irreducible general fibers. A fiber $F$ of $p$ is called {\em singular} (resp. {\em reducible}) 
if $F$ is scheme-theoretically not isomorphic to a general fiber (resp. if $F$ has more than one irreducible components). 
If $F$ is reducible, any irreducible component of $F$ is called a {\em fiber component} of $F$. The fibration $p$ is 
called {\em completely separated} if every reducible fiber is a disjoint union of irreducible fiber components. 
If the morphism $p : Y \to X$ is smooth, then it is completely separated. 

For a $\C$-algebra $B$ and a $\C$-derivation $D$ of $B$, we say that an ideal $I$ of $B$ is {\em $D$-stable} if $D(I) 
\subseteq I$. Similarly, a closed set $V$ of $\Spec B$ is {\em $D$-stable} if the radical defining ideal of $V$ in $B$ 
is $D$-stable. We resume our assumptions and notations. We give a result concerning $D$-stability of the irrdecible 
components of a fiber $F_c$ of $p$. Let $F_c$ be the fiber of $p$ defined by $\xi-c=0$. Let $R=B/(\xi-c)$. Write 
$\xi-c=f_1^{\alpha_1}\cdots f_n^{\alpha_n}$ be the prime decomposition of $\xi-c$ in $B$ with $\forall \alpha_i > 0$. 

\begin{lem}\label{Lemma 2.2}
With the above notations, we have the following assertions.
\begin{enumerate}
\item[(1)]
Let $\gp$ be a minimal prime divisor of $(0)$ in $R$. Then $\gp$ is $D$-stable, and hence $R/\gp$, which is equal 
to $B/f_i$ for some $i$, has the $\C$-derivation induced by $D$. In particular, the nilradical $\sqrt{0}$ is $D$-stable.
\item[(2)]
Every irreducible fiber $F_c$ is reduced. It is also smooth if $D$ is surjective and $B$ is smooth. 
\item[(3)]
Assume that $D$ is surjective. Let $F_c=\Spec B/(\xi-c)$ be an irreducible fiber. Then the $\C$-derivation $D_c$ 
induced on $B/(\xi-c)$ is surjective. 
\end{enumerate}
\end{lem}

\Proof
(1)\ The result is known (see e.g., \cite{F}). But we give a proof for the completeness of the proof.  Let $\gp$ be 
a minimal prime divisor of $(0)$ in $R$. Then $\gp$ is the ideal quotient $(0:a)$ for an element $a \in R$. Namely, 
$\gp=\{z \in R \mid az=0\}$. Note that $(0:a)=(0:a^2)$ since $\gp$ is a prime ideal. Applying $D$ to $az=0$, we have 
$aD(z)+D(a)z=0$, whence $a^2D(z)=0$ and $D(z) \in (0:a^2)=\gp$. If $a \not\in \gp$, then $(0:a)=(0:a^2)$ holds. 
But if $a \in \gp$ then $a^2=0$ and $(0:a^2)=R$. Hence this argument breaks. We need another argument to show 
that $\gp$ is $D$-stable. First, let $(R,M)$ be an Artin local ring over $\C$ with a $\C$-derivation $D$. 
Then the maximal ideal $M$ is $D$-stable. In fact, if $b$ is an element of $M$ with $D(b) \not\in M$. Then $b^n=0$ 
for some $n > 0$ and $D(b)$ is a unit. We have $0=D(b^n)=nb^{n-1}D(b)$, whence $b^{n-1}=0$. Repeating this argument, 
we have $b=0$, which is a contradiction. Now let $P$ be a minimal prime divisor of $(0)$ in a Noetherian ring $R$ 
over $\C$. Consider the local ring $(R_P, PR_P)$. Then it is an Artin local ring defined over $\C$ with a $\C$-derivation 
$D$. By the foregoing argument, $PR_P$ is $D$-stable. Let $b$ be an element of $P$. Then $sD(b) \in P$ for some element 
$s \not\in P$. Then $D(b) \in P$. Hence $P$ is $D$-stable. Let $I$ be a $D$-stable ideal and let $P$ be a minimal prime 
divisor of $R$. We can argue with the residue ring $R/I$ and the induced derivation $\ol{D}$. Then $P/I$ is a minimal 
prime divisor of $(0)$ in $R/I$. Hence it follws that $P$ is $D$-stable. Now we know that any minimal prime ideal $\gp$ 
of $(0)$ is $D$-stable and $D$ induces a $\C$-derivation on the residue ring $R/\gp$. Write 
$\sqrt{0}=\gp_1\cap \gp_2\cap \cdots\cap \gp_n$ with minimal prime divisors $\gp_i$ of $(0)$ in $R$. Since each 
$\gp_i$ is $D$-stable, $\sqrt{0}$ is also $D$-stable. 

(2)\ We show that if $F_c=\Spec R$ is irreducible, it is reduced. In fact, suppose that $\Spec R$ is not reduced. 
Then $\xi-c=f^n$ with $n > 1$. Then $0=D(\xi-c)=nf^{n-1}D(f)$, whence $D(f)=0$ and $f \in \C[\xi]$. Writing $f$ as an 
element in $\C[\xi]$ and substituting it in $\xi-c=f^n$, we obtain an algebraic relation of $\xi$ over $\C$. This is 
a contradiction. Hence $\Spec R$ is reduced. Let $P$ be a singular point of $F_c$ and let $\{x,y\}$ be a local system 
of parameters of $Y=\Spec B$ at the point $P$. Let $f=\xi-c$. Let $f_x, f_y$ be the partial derivatives of $f$ 
with respect to $x, y$ in the local ring $\SO_{Y,P}$ and write $f_x=\ol{f}_xh, f_y=\ol{f}_yh$ with $h=\gcd(f_x,f_y)$ 
in $\SO_{Y,P}$. Since $F_c$ is the curve defined locally by $f=0$, we have $f_xdx+f_ydy=0$ on $F_c$.  
Then $\ol{f}_x(P)=\ol{f}_y(P)=0$ as $P$ is a singular point of $F_c$. On the other hand, we have $f_xD(x)+f_yD(y)=0$ 
in $B$ because $D(f)=D(\xi-c)=0$. Hence $D(x)=\ol{f}_yh_1$ and $D(y)=-\ol{f}_xh_1$ with $h_1 \in \SO_{Y,P}$. 
Now the surjectivity of $D$ implies that there exists an element $t \in B$ such that $D(t)=1$. Then 
$D(t)=t_xD(x)+t_yD(y)=t_x\ol{f}_yh_1-t_y\ol{f}_xh_1$ and $D(t)(P)=0$. This is a contradiction. Hence $F_c$ is smooth.

(3)\ Let $\rho_c : B \to B/(\xi-c)$ be the residue homomorphism. Then $\rho_c\cdot D=D_c\cdot\rho_c$. Hence $D_c$ is 
surjective provided so is $D$.
\QED

Now we can prove the following result which includes Theorem 1 in the introduction.

\begin{thm}\label{Theorem 2.3}
Let $B$ be a factorial affine domain of dimension two over $\C$. Assume that $B^*=\C^*$ and $B$ has a surjective 
$\C$-derivation $D$ with $\Ker D \supsetneqq \C$. Then $B$ is isomorphic to a polynomial ring $\C[x,y]$ and $D$ 
corresponds to the partial derivation $\partial/\partial x$. 
\end{thm}
\Proof
Let $A=\Ker D$. By Lemma \ref{Lemma 2.1}, $A$ is a polynomial ring $\C[\xi]$. Let $Y=\Spec B$, $X=\Spec A$ and 
$p : Y \to X$ be the morphism induced by the natural inclusion $A \hookrightarrow B$. By Lemma \ref{Lemma 2.1}, (2),
$p$ is the fibration. Let $F_c$ be a general smooth fiber and let $R=B/(\xi-c)$. By Lemma \ref{Lemma 2.2}, there 
is a surjective $\C$-derivation on $R$ induced by $D$, which we denote by the same symbol $D$. By Theorem 
\ref{Theorem 1.2}, the fiber $F_c$ is isomorphic to the affine line. Hence $p$ is an $\A^1$-fibration over 
$X=\Spec A\cong \A^1$. In fact, if we choose an element $t \in B$ so that $D(t)=1$, then $R=\C[\:\ol{t}\:]$, where 
$\ol{t}$ is the residue class of $t$ in $R$. Since $B$ is factorial and $B^*=\C^*$, by an algebraic characterization 
of the affine plane \cite{M1}, $Y$ is isomorphic to $\A^2$ and $B$ is isomorphic to $\C[\xi,\tau]$, where we can 
take the generator $\xi$ of $A$ as one of the parameters of $B$. Let $B_0=\C[\xi,t]$, where $D(t)=1$, and let 
$Z=\Spec B_0$. Then the inclusion $B_0 \hookrightarrow B$ induces a dominant $X$-morphism $\pi : Y \to Z$. Since 
$B_0/(\xi-c)=\C[\:\ol{t}\:]=B/(\xi-c)$ with the second equality being proved in the proof of Theorem \ref{Theorem 1.2}, 
the morphism $\pi$ is a bijection. Hence $\pi$ is an isomorphism by Zariski's main theorem and thereby $B=\C[\xi,t]$. 
If we make a change of variables $\xi=y$ and $t=x$, then $D$ is identified with $\partial/\partial x$. 
\QED

The following result gives a characterization of a locally nilpotent derivation on $\C[x,y]$.

\begin{cor}\label{Corollary 2.4}
Let $D$ be a non-trivial reduced $\C$-derivation of $\C[x,y]$. Then $D$ is locally nilpotent if and only if $D$ is surjective 
and $\Ker D \supsetneqq \C$.
\end{cor}

The following remark is also interesting.

\begin{remark}\label{Remark 2.5}{\em
Let $B=\C[x,y]$ be the ring of polynomials in two variables and let $\{f, g\}$ be a Jacobian pair, i.e., the Jacobian 
determinant $\det \partial(f,g)/\partial(x,y)=1$. Let $D$ be a $\C$-derivation 
\[
D=\frac{\partial f}{\partial y}\frac{\partial}{\partial x}-\frac{\partial f}{\partial x}\frac{\partial}{\partial y}.
\]
Then the Jacobian conjecture holds if and only if $D$ is surjective.}
\end{remark}

In fact, $D(f)=0$ and $D(g)=1$. If $D$ is surjective, then Corollary \ref{Corollary 2.4} implies that $D$ is locally 
nilpotent. It is well-known that the polynomial endomorphism of $\A^2$ defined by the pair $\{f,g\}$ is an automorphism 
if $D$ is locally nilpotent. The converse is also clear. For these facts, we refer to \cite[Chapter 4]{F}.

\section{Some observations in general cases}

Let $Y=\Spec B$ be a {\em smooth} affine variety and let $D$ be a surjective $\C$-derivation of $B$. We summarize the results 
obtained in this setting. 

\begin{lem}\label{Lemma 3.1}
The following assertions hold.
\begin{enumerate}
\item[(1)]
Let $\gp$ be a $D$-stable prime ideal of $B$ with $\height(\gp)=\dim X-1$. Then the closed set $V(\gp)$ defined by 
$\gp$ is isomorphic to $\A^1$.
\item[(2)]
Let $\gp_i~(i=1,2)$ be $D$-stable prime ideals of $B$ with $\height(\gp_1)=\dim X-1$. If $\gp_1 \ne \gp_2$ then 
$V(\gp_1)\cap V(\gp_2)=\emptyset$.
\end{enumerate}
\end{lem}
\Proof
(1)\ Let $R=B/\gp$. By the hypothesis, $\dim R=1$. Since $\gp$ is $D$-stable, $D$ induces a surjective $\C$-derivation 
$\ol{D}$ on $R$. By Proposition \ref{Proposition 1.7} and Theorem \ref{Theorem 1.2}, it follows that $V(\gp)\cong \A^1$.

(2)\ Suppose that $V(\gp_1)\cap V(\gp_2) \ne \emptyset$. Let $P \in V(\gp_1)\cap V(\gp_2)$. Then $P$ corresponds to 
a maximal ideal $\gm$ of $B$ such that $\gm \supset \gp_1+\gp_2$. Since $\height(\gp_1)=\dim X-1$, we know 
that $\gm$ is a prime divisor of $\gp_1+\gp_2$. Since $\gp_1+\gp_2$ is $D$-stable, $\gm$ is also $D$-stable (cf. 
the proof of Lemma \ref{Lemma 2.2}). Then $D$ induces a non-trivial $\C$-derivation on $B/\gm$ because $D(t)=1$. 
This is a contradiction because $D$ is trivial on $\C$.
\QED

\begin{cor}\label{Corollary 3.2}
With the notations and assumptions in Lemma \ref{Lemma 3.1}, assume further that $A:=\Ker D$ is finitely generated over 
$\C$ and $\dim A=\dim B-1$. Let $X=\Spec A$ and let $p : Y \to X$ be the morphism associated to the inclusion 
$A \hookrightarrow B$. Then $p$ is an $\A^1$-fibration. Furthermore, every one-dimensional fiber component of $p$ 
is the affine line and disjoint from other fiber components.
\end{cor}
\Proof
Since $Q(A)$ is algebraically closed in $Q(B)$ (cf. the proof of Lemma \ref{Lemma 2.1}), $p$ is an $\A^1$-fibration. 
In fact, for a maximal ideal $\gm$ of $A$, the ideal $\gm B$ is $D$-stable, and hence every minimal prime divisor of 
$\gm B$ is $D$-stable. Consider the fiber $F$ of $p$ over the closed point of $X$ corresponding to $\gm$. If one fiber 
component $F_1$ of $F$ has dimension one, then it is the affine line by Lemma \ref{Lemma 3.1}, (1) and disjoint from 
the other fiber components of $F$ by Lemma \ref{Lemma 3.1},(2). In particular, any general closed fiber is the affine line 
since $\gm B$ is then a prime ideal. By restricting $X$ to an open set $D(f)$ of $X$ with $f \in A$ and replacing $X$ 
by $\Spec B[f^{-1}]$, we may assume that $p$ is faithfully flat, $X$ is smooth and every fiber of $p$ is irreducible 
and reduced. Then $p$ is an $\A^1$-bundle by \cite[Theorem 2]{KM}.
\QED

In a trial of reproving Cerveau's \lq\lq theorem\rq\rq~in the case $\Ker D=\C$, we can prove the following result.

\begin{thm}\label{Theorem 3.3}
Let $R$ be a normal affine domain of dimension one over $\C$ with a $\C$-derivation $D$. Let $t$ be an indeterminate 
and let $B=R[t]$. Extend the derivation $D$ to a $\C$-derivation $\wt{D}$ on $B$ by setting $\wt{D}(t)=1$. Assume 
that the following conditions are satisfied:
\begin{enumerate}
\item[(1)]
$\wt{D}$ is surjective.
\item[(2)]
$\Ker\wt{D}=\C$.
\item[(3)]
There exists a maximal ideal $\gm$ of $R$ such that $D(\gm) \subseteq \gm$.
\end{enumerate}
Then $R$ is a polynomial ring $\C[x]$ and $D=cx(d/dx)$ after a suitable choice of $x$ and with $c \in \C^*$. 
Hence $B=\C[x,t]$ and $\wt{D}=cx(\partial/\partial x)+(\partial/\partial t)$.
\end{thm}
\Proof
We first prove the following
\svskip

\noindent
{\bf Claim 1.}\ {\em Let $\varphi : R \to \C$ be the natural residue homomorphism of $R$ to $\C=R/\gm$. Then $D$ 
induces a bijective $\C$-linear mapping $\gm \stackrel{\sim}{\to} \gm$.}
\svskip

For an element $a \in R$, there exists an element $h(t) \in B$ such that $\wt{D}(h(t))=a$. Write 
\[
h(t)=a_0+a_1t+\cdots+a_nt^n,~~a_i \in R.
\]
Then we have 
\[
a=D(a_0)+a_1,~~D(a_{i-1})+ia_i=0~(1<i\le n),~~D(a_n)=0.
\]
Hence $a_n \in \Ker D$. By the condition (2), write $a_n=-c_n \in \C$. Then $D(a_{n-1})+na_n=0$ implies that 
$\wt{D}(a_{n-1}-nc_nt)=0$, whence $a_{n-1}=nc_nt+c_{n-1}$. Since $t$ is an indeterminate over $R$, it follows that 
$c_n=0$. Similarly, $D(a_{n-2})+(n-1)a_{n-1}=0$ implies that $c_{n-1}=0$. This argument applies until the relation 
$D(a_1)+2a_2=0$, and it follows that $a_2=\cdots=a_n=0$ and $a_1 \in \C$. Hence we have $a=D(a_0)+a_1$. The 
element $a_1$ is uniquely determined by $a$ and $D(a_0) \in \gm$. In fact, if we have another expression 
$a=D(a'_0)+a'_1$ with $a'_0 \in R$ and $a'_1 \in \C$, then $D(a_0-a'_0)=a'_1-a_1$. Hence $a_0-a'_0=(a'_1-a_1)t+c$ 
with $c \in \C$. This implies $a'_1=a_1$. Write $\rho(a)=a_1$. On the other hand, write $a_0=a'+c$ with 
$a' \in \gm$ and $c \in \C$. Then $D(a_0)=D(a') \in \gm$ by the condition (3). Hence the decomposition 
$a=D(a_0)+a_1$ is given by the $\C$-module decomposition $R=\gm\oplus \C$. Define a mapping $\rho : R \to \C$ 
by $\rho(a)=a_1$. Since $\gm$ is an ideal, $\rho$ is a ring homomorphism and $D$ induces a $\C$-linear mapping 
$D : \gm \to \gm$. It is injective since $\Ker D \cap \gm =(0)$ by the condition (2). Furthermore, if $a \in \gm$, 
the above decomposition $a=D(a_0)+a_1$ gives $a_1=0$. Namely, $a=D(a_0)$. Thus $D$ is surjective on $\gm$. 
We then prove the next
\svskip

\noindent
{\bf Claim 2.}\ {\em Let $\delta : \Omega_{R/\C} \to R$ be the $R$-module homomorphism corresponding to 
the derivation $D$. Then $\delta$ induces an isomorphism between $\Omega_{R/\C}$ and $\gm$. Furthermore, every 
$1$-form $\omega \in \Omega_{R/\C}$ is exact.}
\svskip

In fact, since $d(a)=d(a+c)$ for every $c \in \C$, it follows that $\Omega_{R/\C}$ is generated by $\{da \mid 
a \in \gm\}$ as the $R$-module. Since $\delta(da)=D(a)$, it follows that $\im \delta$ is the $R$-ideal generated 
by $\{D(a) \mid a \in \gm\}$. By Claim 1, we know that $\im \delta=\gm$. Since $\delta$ is clearly injective, 
$\delta$ induces an $R$-isomorphism between $\Omega_{R/\C}$ and $\gm$. Let $\omega \in \Omega_{R/\C}$. Then 
$\delta(\omega)=D(a)$ with $a \in \gm$. Hence $\delta(\omega-da)=0$, and $\omega=da$. 

Now, by the proof of Theorem \ref{Theorem 1.2}, $R$ is a polynomial ring $\C[x]$. Write $\gm=(x)$ after a suitable 
change of the variable $x$. Then $D(x)=xf(x)$ with $f(x) \in \C[x]$. Since $D : \gm \to \gm$ is bijective, there 
exists $h(x) \in \gm$ such that $D(h(x))=x$. Then $h'(x)f(x)=1$. Hence $f(x)=c \in \C^*$ and $D=cx(d/dx)$.
Then it is clear that $\wt{D}=cx(\partial/\partial x)+(\partial/\partial t)$.
\QED

If we drop the condition (3) in the assumptions of Theorem \ref{Theorem 3.3}, $R$ is no longer a polynomial ring 
as shown by the following example. 

\begin{example}\label{Example 3.4}
Let $R$ be a Laurent polynomial ring $\C[x,x^{-1}]$ and let $D=x(d/dx)$. Let $B=R[t]$ be a polynomial ring in a 
variable $t$ over $R$. Extend $D$ to a $\C$-derivation $\wt{D}=x(\partial/\partial x)+(\partial/\partial t)$. 
Then $\wt{D}$ is surjective and $\Ker \wt{D}=\C$, but there is no maximal ideal $\gm$ of $R$ such that 
$D(\gm) \subseteq \gm$.
\end{example}

If the condition (3) above is not satisfied, we show that we are essentially in the situation of Example 
\ref{Example 3.4}.

\begin{thm}\label{Theorem 3.5}
Let $R$ be a normal affine domain of dimension one over $\C$ with a $\C$-derivation $D$. Let $t$ be an indeterminate 
and let $B=R[t]$. Extend the derivation $D$ to a $\C$-derivation $\wt{D}$ on $B$ by setting $\wt{D}(t)=1$. Assume 
that the following conditions are satisfied:
\begin{enumerate}
\item[(1)]
$\wt{D}$ is surjective.
\item[(2)]
$\Ker\wt{D}=\C$.
\item[(3)]
There are no maximal ideals $\gm$ of $R$ such that $D(\gm) \subseteq \gm$.
\end{enumerate}
Then $R$ is a Laurant polynomial ring $\C[x,x^{-1}]$ and $D=cx(d/dx)$ after a suitable choice of $x$ and with $c \in \C^*$. 
Hence $B=\C[x,x^{-1},t]$ and $\wt{D}=cx(\partial/\partial x)+(\partial/\partial t)$.
\end{thm}
\Proof
Since $\wt{D}$ is surjective and $\Ker \wt{D}=\C$, the proof of Claim 1 above shows that $R=D(R)+\C$, where $D(R)$ is 
the image $\im D$. Since any maximal ideal of $R$ is not $D$-stable by the hypothesis, the ideal of $R$ generated by 
$\{D(a) \mid a \in R\}$ is the unit ideal. Hence the $R$-homomorphism $\delta : \Omega_{R/\C} \to R$ corresponding to 
the derivation $D$ is surjective. Since $\Omega_{R/\C}$ is a locally free $R$-module of rank $1$, it follows that 
$\delta$ is an isomorphism. Since the $\C$-vector space $D(R)$ has codimension one in $R$, the $\C$-vector space 
$d(R)=\{da \mid a \in R\}$ in $\Omega_{R/\C}$ has also codimension one. By Grothendieck's algebraic de Rham theorem 
\cite{G}, $H^1(X;\C)=\C$, where $X=\Spec R$. Hence $b_1(X)=1$. Then $X \cong \A^1_*=\Spec \C[x,x^{-1}]$, and the 
derivation $D$ is of the form $D=cx^m(d/dx)$, where $c \in \C^*$ and $m \in \Z$, for the derivation has no zeros. If 
$m \le 0$, then $D((c^{-1}/(-m+1))x^{-m+1})=\wt{D}(t)=1$ and $(c^{-1}/(-m+1))x^{-m+1}-t \in \C$, which contradicts 
the hypothesis that $t$ is an indeterminate over $R$. If $m \ge 2$, then $x \not\in D(R)$, which is also 
a contradiction to $R=D(R)+\C$. Hence $m=1$, and $\wt{D}=cx(\partial/\partial x)+(\partial/\partial t)$.
\QED 

Let $B$ be an affine normal domain of dimension two over $\C$ and let $D$ be a $\C$-derivation on $B$. An {\em integral 
curve} with respect to $D$ is an irreducible curve $C$ of $\Spec B$ defined by a prime ideal $\gp$ of height one such 
that $D(\gp) \subseteq \gp$. A nonzero element $f$ of $B$ is an {\em integral element} with respect to $D$ if $D(f)=fh$ 
for some $h \in B$. We say that $D$ has no {\em non-constant multiplicative characters} if $D(f)=fh$ with 
$f, h \in B\setminus \{0\}$ implies $h \in \C^*$. 

\begin{lem}\label{Lemma 3.6}
Let $B$ be the polynomial ring $\C[x,y]$ and let $D$ be a surjective $\C$-derivation on $B$. Then the following 
assertions hold. 
\begin{enumerate}
\item[(1)]
Assume that an integral element $f$ of $B$ exists. Then $f$ is written as a polynomial $f=F(x_2) \in \C[x_2]$, where 
$x_2$ is a coordinate of $B$ and $x_2$ is an integral element. 
\item[(2)]
Assume that $\Ker D=\C$ and that $D$ has no non-constant multiplicative characters. Then there exists at most one 
integral curve on $\Spec B$. 
\item[(3)]
With the same assumptions as in {\rm (2)} above, $\Ker D_K=\C$, where $K$ is the quotient field of $B$ and $D_K$ is 
the extension of $D$ to $K$.
\end{enumerate}
\end{lem}
\Proof
(1)~Let $f=p_1^{\alpha_1}\cdots p_n^{\alpha_n}$ be the irreducible decomposition of $f$. Then $p_1, \ldots, p_n$ 
are integral elements by Lemma \ref{Lemma 2.2}. Since $\Spec B/(p_i)$ $\cong \A^1$, it follows that $p_i$ is 
a coordinate of $B$. We write $p_1=x_2$ with $B=\C[x_1,x_2]$. For $i \ne 1$, $V(p_1)\cap V(p_i)=\emptyset$. 
Otherwise, there exists a $D$-stable maximal ideal of $B$, which is impossible. Write $p_i=x_2g(x_1,x_2)+\varphi(x_1)$.
Then $\varphi(x_1)=c \in \C^*$. Since $p_i$ is also a coordinate, $p_i-c$ is irreducible and hence $g(x_1,x_2) \in \C^*$.
So, $f$ is a polynomial in $x_2$.

(2)~Let $\gp_i=(u_i)$ define an integral curve for $i=1,2$. By (1) above, we may assume that $u_1=x_2$ is a coordinate
and $u_2=x_2+c$ with $c \in \C^*$. Since $D(x_2)$ is divisible by $x_2$, either $D(x_2)=0$ or $D(x_2)=x_2h $ with 
$h \in B\setminus \{0\}$. If $D(x_2)=0$ then $x_2 \in \Ker D=\C$, which is a contradiction. In the latter case, 
since $D$ has no non-constant multiplicative characters, it follows that $D(x_2)=\alpha x_2$ with $\alpha \in \C^*$. 
Replacing $D$ by $\alpha^{-1}D$, we may assume that $D(x_2)=x_2$. Since $D(x_2)=D(x_2+c)=\beta (x_2+c)$ with 
$\beta \in \C^*$, we have $x_2=\beta(x_2+c)$, which is a contradiction. 

(3)~Let $\xi \in \Ker D_K$. Write $\xi=f/g$, where $f, g \in B$ and $\gcd(f,g)=1$. If $g \in \C^*$, then $\xi \in 
\Ker D$ and hence $\xi \in \C^*$. Suppose that $g \not\in \C^*$. Then $D(f) \in (f)$ and $D(g) \in (g)$. Every 
irreducible component of $f, g$ is also $D$-stable. Hence there exist irreducible $D$-stable elements $p, q$ 
of $B$ such that $p \mid f$ and $q \mid g$. Since there is at most one integral curve on $\Spec B$, it follows 
from the foregoing assertion that $q=cp$ with $c \in \C^*$. This is a contradiction.  
\QED

Given a $\C$-derivation $D$ on an affine domain $B$, the existence of an integral element is not clear. 
When $B$ is a polynomial ring $\C[x,y]$, $D$ being surjective implies the existence of an integral element. We 
show this result after Cerveau's argument \cite{cerveau}, where an essential ingredient is a theorem of 
Dimca-Saito \cite{DS}. We denote by $\C[x_1,\ldots, x_n]$ the polynomial ring $B$ in dimension $n$, by $\Omega^i$ 
the free $B$-module of differential $i$-forms $(0 \le i \le n)$ and by $\Omega^\bullet$ the differential 
complex $\{\Omega^i~(0 \le i \le n); d : \Omega^i \to \Omega^{i+1}\}$, where $d$ is defined by the exterior 
differentiation
\[
d(fdx_{k_1}\wedge\cdots\wedge dx_{k_i})=\sum_{j=1}^n\frac{\partial f}{\partial x_j}dx_j\wedge dx_{k_1}\wedge 
\cdots\wedge dx_{k_i}.
\]
Given $f \in B$, define another differentiation $D_f$ of $\Omega^\bullet$ by 
\[
D_f(\omega)=d\omega +df\wedge \omega, \quad \omega \in \Omega^i.
\]
Then we have a crucial result of Dimca-Saito \cite{DS}.

\begin{lem}\label{Lemma 3.7}
Let $f \in B$ and denote by $\varphi_f : \A^n \to \A^1$ the morphism induced by the inclusion $\C[f] 
\hookrightarrow B$. Let $F$ be a general fiber of the morphism $\varphi_f$. Then there is an isomorphism 
\[
H^{i+1}(\Omega^\bullet, D_f) \cong \wt{H}^i(F; \C) \quad \text{for every $i$},
\]
where $\wt{H}^i(F; \C)$ denotes the reduced cohomology. 
\end{lem}

Now we consider the case $n=2$, which is due to Cerveau \cite{cerveau} and stated in the introduction as Theorem 2. 

\begin{thm}\label{Theorem 3.8}
Let $D=f_1(\partial/\partial x_1)+f_2(\partial/\partial x_2)$ be a surjective derivation on $B=\C[x_1,x_2]$. 
Assume that $\Ker\delta=\C$. Then there exists an integral element with repsect to $D$.
\end{thm}
\Proof
Set $\mathrm{vol}=dx_1\wedge dx_2$ (the volume form) and set 
\[
\omega_D=i_D(\mathrm{vol})=-f_2dx_1+f_1dx_2.
\]
Then we have 
\[
d\omega_D=\left(\frac{\partial f_1}{\partial x_1}+\frac{\partial f_2}{\partial x_2}\right) dx_1\wedge dx_2.
\]
Since $D$ is surjective, there exists an element $g \in B$ such that 
\[
D(g)=f_1\frac{\partial g}{\partial x_1}+f_2\frac{\partial g}{\partial x_2}=
-\left(\frac{\partial f_1}{\partial x_1}+\frac{\partial f_2}{\partial x_2}\right).
\]
Then we have 
\[
d\omega_D+dg\wedge \omega_D=0,\quad \mathrm{i.e.}, D_g(\omega_D)=0.
\]
Suppose that the curve $g=c$ is irreducible for a general $c \in \C$. Then a general fiber $G$ of the morphism 
$\varphi_g : \A^2 \to \A^1$ is irreducible. Hence $\wt{H}^0(G; \C)=0$. By Lemma \ref{Lemma 3.7}, there exists 
an element $h \in B$ such that $D_g(h)=\omega_D$. Since $D_g(h)=d(h)+hdg$, we have 
\[
\omega_D=d(h)+hdg.
\]
Note that the derivation $D$ is given as $D=\delta\cdot d$ for a $B$-module homomorphism $\delta : \Omega^1 \to 
B$. In fact, $\delta(dx_1)=f_1$ and $\delta(dx_2)=f_2$. Hence we have $\delta(\omega_D)=-f_2f_1+f_1f_2=0$. 
Thence we have 
\[
0=\delta(\omega_D)=\delta(d(h)+hdg)=D(h)+hD(g).
\]
Thus the element $h$ is an integral element. Finally, the assumption that $g=c$ is irreducible for general $c 
\in \C$ (the generic irreducibility of $g$) is guaranteed by Cerveau \cite[Section 6]{cerveau}. 
\QED

By Lemma \ref{Lemma 3.6}, (1), the integral element $h \in B$ is an element of $\C[x_2]$, where $x_2$ is one of the 
coordinates. Since $\omega_D$ has no divisorial part, the polynomial $h(x_2)$ has only simple factors. Cerveau 
\cite[Lemme 5.3]{cerveau} asserts that $h(x_2)$ is a linear polynomial. Thus we may assume that $\omega_D= 
c\omega'$ after a change of coodinates $x_2 \mapsto x_2-\alpha$ with $\alpha \in \C$, where $\omega'=dx_2+x_2dg$. 

\begin{lem}\label{Lemma 3.9}
With the assumptions and notations in the proof of Theorem \ref{Theorem 3.8}, the following assertions hold.
\begin{enumerate}
\item[(1)]
For every $\eta \in \Omega^1$, there exist $f \in B$ and $\eta_D \in \Ker\delta$ such that $\eta=\eta_D+df$ 
{\em (see \cite[Proposition 1.2]{cerveau}).} Furthermore, $\eta\wedge \omega_D=\delta(\eta)dx_1\wedge dx_2$.
\item[(2)]
$\Ker \delta=B\omega_D$.   
\item[(3)]
Let $\eta=h_1dx_1+h_2dx_2$ be an element of $\Ker\delta$. If $\eta$ has no divisorial part, then $\eta=c\omega_D$ 
with $c \in \C^*$.
\item[(4)]
Let $\eta=f\omega_D$ be an element of $\Ker\delta$. Then we have
\[
d\eta+dg\wedge \eta=df\wedge \omega_D=D(f)dx_1\wedge dx_2.
\]
Hence $d\eta+dg\wedge \eta=0$ if and only if $f \in \C$.
\end{enumerate}
\end{lem}
\Proof
(1)~~The module $\Ker\delta$ is a free $B$-submodule of rank one of $\Omega^1$ which is a direct summand of $\Omega^1$. 
This follows from the following exact sequence
\[
0 \to \Ker\delta \to \Omega^1 \st{\delta} B \to 0.
\]
Take $\eta=h_1dx_1+h_2dx_2$ from $\Omega^1$. Since $D$ is surjective, there exists $f \in B$ such that $D(f)=h_1f_1
+h_2f_2$. Set $\eta_D=\eta-df$. Then $\eta_D \in \Ker\delta$. The rest is easy to show.

(2)~~Let $\eta=h_1dx_1+h_2dx_2$ be an element of $\Ker\delta$. Then $a\eta=b\omega_D$ with $a, b \in B$, $a \ne 0$ and 
$\gcd(a,b)=1$. Then $ah_1=-bf_2$ and $ah_2=bf_1$. Since $\gcd(a,b)=1$, we have $b \mid h_1$ and $b \mid h_2$. Writing 
$h_1=bc$ and $h_2=bd$ with $c, d \in B$, we have $f_1=ad$ and $f_2=-ac$. Since $\gcd(f_1,f_2)=1$ as $D$ is 
surjective, it follows that $a \in \C^*$. Hence $h_1=-ba^{-1}f_2$ and $h_2=ba^{-1}f_1$. Namely, $\eta=a^{-1}b\omega_D$.

(3)~~Suppose that $\eta$ has no divisorial part, i.e., $\gcd(h_1,h_2)=1$. Then $b \in \C^*$ as well. 
Hence $a^{-1}b \in \C^*$. 

(4)~~Since we have $d\eta+dg\wedge\eta=df\wedge\omega_D+f(d\omega_D+dg\wedge\omega_D)=df\wedge\omega_D=D(f)dx_1\wedge dx_2$,
it is easy to obtain the conclusion.
\QED

Replacing $D$ by $cD$ with $c \in \C^*$ if necessary, we may assume that $\omega_D=\omega'=dx_2+x_2dg$. 
It is shown in \cite[Remarque 5.6, Lemmes 5.4 and 5.5]{cerveau} that $g$ is of the form $W=\lambda x_1+\varphi(x_2)$ 
with $\lambda \in \C^*$ and $\varphi(x_2) \in \C[x_2]$. But there is an elementary (computational) mistake in 
the proof of Lemma 5.5. Thus the assertion of Cerveau is yet to be proved. 

Nowithstanding, we can show the following result.

\begin{lem}\label{Lemma 3.10}
With the above notations, especially with the notations in the proof of Theorem \ref{Theorem 3.8}, the fibration 
$\varphi_g : \A^2 \to \A^1$ is an $\A^1$-fibration. Hence $g$ is a ring generator of $B$.
\end{lem}
\Proof
We use the result of Dimca-Saito and prove that the $B$-module homomorphism $D_g : \Omega^1_{B/\C} \to \Omega^2_{B/\C}$ 
is surjective. Since $\Omega^2=Bdx_1\wedge dx_2$, it suffices to show that $\varphi dx_1\wedge dx_2 \in \im D_g$ 
for every $\varphi \in B$. Since $D$ is surjective, we can choose $f \in B$ so that $\varphi=D(f)$. Set $\eta=f\omega_D$. 
By Lemma \ref{Lemma 3.9}, (4), we have $D_g(\eta)= d\eta+dg\wedge \eta= D(f)dx_1\wedge dx_2=\varphi dx_1\wedge dx_2$. 
Hence $\Omega^2=\im D_g$. Let $G$ be a general fiber of $\varphi_g$. By the foregoing argument, $G$ is a smooth curve 
and $H^1(G;\C)=0$ by Lemma \ref{Lemma 3.7}. Hence $G$ is isomorphic to $\A^1$, and $\varphi_g$ is an $\A^1$-fibration.
\QED

We summarize various obstacles as questions.

\begin{question}\label{Question 3.11}{\em
Let $W \in B$. Is the ideal of $B$ generated by $x_2\frac{\partial W}{\partial x_1}$ and 
$1+x_2\frac{\partial W}{\partial x_2}$ a unit ideal ? If necessary, we assume the condition that $D(x_2)=-x_2D(W)$. }
\end{question}
If the answer is positive, then $dx_2+x_2dW=c\omega_D$ with $c \in \C^*$ and $W-g$ is a constant.

\begin{question}\label{Question 3.12}{\em
Let $D$ be a surjective $\C$-derivation on $\C[x_1,x_2]$. Does $D$ then have no non-constant multiplicative characters ?}
\end{question}
If the assertion of Cerveau is verified, it gives a positive answer to the second question. 
\svskip

To finish this article, we give several examples of $\C$-derivations on $B=\C[x,y]$ which have no integral elements. We change 
here again the notations of variables from $x_1, x_2$ to $x, y$. This is to accord with the notations in the original sources 
which we refer to. A. Nowicki and other people \cite{A, A2, OA} gave such examples of a derivation 
$D=(\partial/\partial x)+f(x,y)(\partial/\partial y)$ on $B$. In fact, $D$ has no integral curves if $f(x,y)$ is one 
of the following polynomials
\begin{enumerate}
\item[1]
$f(x,y)=xy+1$,
\item[2]
$f(x,y)=(x^2+x)y+x^2$,
\item[3]
$f(x,y)=y^s+cx$ with $ s \ge 2$ and $c \in \C^*$.
\end{enumerate}

The following example is due to L. Makar-Limanov \cite{L} and it seems to be new in construction. So, we include the computation 
which Makar-Limanov kindly communicated to us. 

\begin{example}\label{Example 3.13}
Let $D=(x^2+\lambda y)(\partial/\partial x)+(y^2+\mu x)(\partial/\partial y)$ on $\C[x,y]$, where $\lambda, \mu \in 
\C^*$ with $\lambda^2 \ne \mu^2$. Then $D$ has no integral elements. 
\end{example}
\Proof
Assume that $D(f)=gf$. Then $\deg_{x,y}(g) < 2$, where $\deg_{x,y}$ here signifies the total degree. Hence we may assume that 
$g=\alpha x+\beta y+ \gamma$. Write $f=x^af_0+r$, where $\deg_x(r) < a$ and $f_0 \in \C[y]$. Then 
\[
D(f)=(af_0+\mu f'_0)x^{a+1}+\cdots =\alpha f_0x^{a+1}+\cdots
\]
and hence $af_0+\mu f'_0=\alpha f_0$. Hence $\alpha=a$ and $f'_0=0$. Similarly, $\beta=b=\deg_y(f)$ and $f=c'y^b+r'$ with 
$c' \in \C$ and $\deg_y(r') < b$.

Consider now the leading form $\ol{f}$ of $f$ relative to the total degree. Write $\ol{f}=x^iy^jh$, where 
$h(x,0)h(0,y) \ne 0$. We have $\ol{D}(\ol{f})=(ax+by)\ol{f}$, where $\ol{D}(x)=x^2$ and $\ol{D}(y)=y^2$. Since 
$\ol{D}(x^iy^j)=(ix+jy)x^iy^j$, we have 
\[
\ol{D}(h)=\{(a-i)x+(b-j)y\}h.
\]
Since $h=h_0x^\ell+\cdots+h_\ell y^\ell$ with $h_0h_\ell \ne 0$, we have $a-i=\ell, b-j=\ell$ and $\ol{D}(h)=
\ell(x+y)h$. It is easy to see that all coefficients of $h$ are defined uniquely by $h_0$ and that $h_0(x-y)^\ell=h$. Hence 
$\ol{f}=h_0x^iy^j(x-y)^\ell$ and $a=\deg_x(f)=i+\ell, b=\deg_y(f)=j+\ell$. Since we know that the leading monomials of 
$f$ relative to $\deg_x$ and $\deg_y$ are respectively $cx^a$ and $c'y^b$ with $c, c'\in \C$, we can conclude that 
\[
i=j=0,~~a=b=\ell,~~\ol{f}=(x-y)^a,~~D(f)=\{a(x+y)+\gamma\}f\ .
\]
Introduce new variables $u=x-y$ and $v=x+y$. Then 
\begin{eqnarray*}
D(u)&=& uv+\frac{\lambda-\mu}{2}v-\frac{\lambda+\mu}{2}u, \\
D(v)&=& \frac{1}{2}\{u^2+v^2+(\lambda+\mu)v+(\mu-\lambda)u\}, \\
D(f)&=& (av+\gamma)f,~~\ol{f}=u^a,~~\deg(f)=a\ .
\end{eqnarray*}
Write $f=v^d\phi+\rho$, where $\deg_v(\rho)< d$. Then 
\[
\frac{d}{2}\phi+\left(u+\frac{\lambda-\mu}{2}\right)\phi'=a\phi
\]
as the coefficients with $v^{d+1}$ in the left and right sides of $$D(f)=(av+\gamma)f.$$ Therefore we have 
\[
\phi=c\left(u+\frac{\lambda-\mu}{2}\right)^{a-\frac{d}{2}}\quad \text{and}\quad 
f=cv^d\left(u+\frac{\lambda-\mu}{2}\right)^{a-\frac{d}{2}}+\rho\ .
\]
Hence $a=\deg(f) \ge d+a-\frac{d}{2}$ and $d = 0$. It means that 
\[
f=f(u)\quad\text{and}\quad D(f)=f'\left(uv+\frac{\lambda-\mu}{2}v-\frac{\lambda+\mu}{2}u\right)=(av+\gamma)f\ .
\]
So, we have 
\begin{eqnarray*}
\lefteqn{\frac{2uv+(\lambda-\mu)v-(\lambda+\mu)u}{av+\gamma}} \\
&&=\frac{2u+\lambda-\mu}{a}-\frac{(2u+\lambda-\mu)\gamma+a(\lambda+\mu)u}{a(av+\gamma)} \in \C(u)
\end{eqnarray*}
which is possible only if $\gamma(\lambda-\mu)=2\gamma+a(\lambda+\mu)=0$. Since we assume that $\lambda^2 \ne \mu^2$, 
we have $\gamma=a=0$. Namely $f \in \C$. 
\QED

\end{document}